\documentclass[11pt,a4paper]{article}

\usepackage{inputenc}
\usepackage{amsmath}
\usepackage{bm}
\usepackage{bbold}
\usepackage{amsthm}

\usepackage{hyperref}

\title{An Algebraic Approach to Project Schedule Development Under Precedence Constraints\thanks{International Journal of Applied Mathematics and Informatics, 2012. Vol. 6, no. 2, pp. 92-100.}
}

\author{Nikolai Krivulin\thanks{Faculty of Mathematics and Mechanics, St.~Petersburg State University, 28 Universitetsky Ave., St.~Petersburg, 198504, Russia, 
nkk@math.spbu.ru}
}

\date{}

\newtheorem{theorem}{Theorem}
\newtheorem{lemma}{Lemma}

\begin{document}

\maketitle

\begin{abstract}
An approach to schedule development in project management is developed within the framework of idempotent algebra. The approach offers a way to represent precedence relationships among activities in projects as linear vector equations in terms of an idempotent semiring. As a result, many issues in project scheduling reduce to solving computational problems in the idempotent algebra setting, including linear equations and eigenvalue-eigenvector problems. The solutions to the problems are given in a compact vector form that provides the basis for the development of efficient computation procedures and related software applications.
\\

\end{abstract}

\section{Introduction}

The problem of scheduling a large-scale set of activities is a key issue in project management \cite{PMI08Aguide,Neumann03Project}. There is a variety of project scheduling techniques developed to handle different aspects of the problem. The techniques range from the classical Critical Path Method and the Program Evaluation and Review Technique marked the beginning of the active research in the area in 1950s, to more recent approaches including methods and techniques of idempotent algebra (see, e.g., \cite{Cuninghame79Minimax,Zimmermann81Linear,Baccelli92Synchronization,Cuninghame94Minimax,Heidergott05Maxplus,Krivulin09Idempotent,Krivulin11Methods} and references therein).

We describe a new computational approach to project scheduling problems, which is based on implementation and further development of models and methods of idempotent algebra in \cite{Krivulin06Solution,Krivulin06Eigenvalues,Krivulin09Idempotent,Krivulin11Methods}. The approach offers a useful way to represent different types of precedence relationships among activities in a project as linear vector equations written in terms of an idempotent semiring. As a result, many issues in project scheduling reduce to solving computational problems in the idempotent algebra setting, including linear equations and eigenvalue-eigenvector problems. We give solutions to the problems in a compact vector form that provides a basis for the development of efficient computation algorithms and related software applications. 

The paper extends previous research presented in \cite{Krivulin12Algebraic} and is organized as follows. We start with a brief introduction to idempotent algebra, that provides main definitions and notation, and then outlines basic results underlying subsequent applications. Furthermore, examples of actual problems in project scheduling are considered. We show how to formulate the problems in terms of idempotent algebra, and present related algebraic solutions. To illustrate the application of the results, numerical examples are given.

\section{Definitions and Notation}

We start with a brief introduction to idempotent algebra based on \cite{Krivulin06Solution,Krivulin06Eigenvalues,Krivulin09Idempotent,Krivulin11Methods}. Further details on the topic can be found in \cite{Cuninghame79Minimax,Zimmermann81Linear,Baccelli92Synchronization,Kolokoltsov97Idempotent,Cuninghame94Minimax,Golan03Semirings,Heidergott05Maxplus,Butkovic10Maxlinear,Litvinov12Idempotent}.

\subsection{Idempotent Semifield}

Consider a set $\mathbb{X}$ that is equipped with two operations $\oplus$ and $\otimes$ called addition and multiplication, and that has neutral elements $\mathbb{0}$ and $\mathbb{1}$ called zero and identity. We suppose that $\langle\mathbb{X},\mathbb{0},\mathbb{1},\oplus,\otimes\rangle$ is a commutative semiring, where addition is idempotent and multiplication is invertible. Since the nonzero elements in $\mathbb{X}$ form a group under multiplication, this semiring is often referred to as the idempotent semifield.

The idempotent property is given by the equality
$$
x\oplus x=x
$$
that is true for all $x\in\mathbb{X}$. Let $\mathbb{X}_{+}=\mathbb{X}\setminus\{\mathbb{0}\}$. For each $x\in\mathbb{X}_{+}$, there exists its inverse $x^{-1}$ such that $x\otimes x^{-1}=\mathbb{1}$.

The power notation is defined as usual. For any $x\in\mathbb{X}_{+}$ and integer $p>0$, we have $x^{0}=\mathbb{1}$, $\mathbb{0}^{p}=\mathbb{0}$, and
$$
x^{p}
=
x^{p-1}\otimes x
=
x\otimes x^{p-1},
\qquad
x^{-p}=(x^{-1})^{p}.
$$

It is assumed that in the semiring, the integer power can naturally be extended to the case of rational exponents.

In what follows, the multiplication sign $\otimes$ is omitted as is usual in conventional algebra. The power notation is thought of as defined in terms of idempotent algebra. However, when writing exponents, we routinely use ordinary arithmetic operations.

Since the addition is idempotent, it induces a partial order $\leq$ on $\mathbb{X}$ according to the rule: $x\leq y$ if and only if $x\oplus y=y$.  With this definition, it is easy to verify that both addition and multiplication are isotonic, and that
$$
x
\leq
x\oplus y,
\qquad
y
\leq
x\oplus y.
$$

The relation symbols are understood below in the sense of this partial order. Note that according to the order, we have $x\geq\mathbb{0}$ for any $x\in\mathbb{X}$.

As an example of the semirings under study, one can consider the idempotent semifield of real numbers 
$$
\mathbb{R}_{\max,+}
=
\langle\mathbb{R}\cup\{-\infty\},-\infty,0,\max,+\rangle.
$$

The semiring has neutral elements $\mathbb{0}=-\infty$ and $\mathbb{1}=0$. For each $x\in\mathbb{R}$, there exists its inverse $x^{-1}$, which is equal to $-x$ in ordinary arithmetics. For any $x,y\in\mathbb{R}$, the power $x^{y}$ is equivalent to the arithmetic product $xy$. The partial order coincides with the natural linear order on $\mathbb{R}$.

We use this semiring as the basis for the development of algebraic solutions to project scheduling problems in the subsequent sections.

\subsection{Vector and Matrix Algebra}

Vector and matrix operations are routinely introduced on the basis of the scalar operations. Consider a Cartesian product $\mathbb{X}^{n}$ with its elements represented as column vectors. For any two vectors $\bm{a}=(a_{i})$ and $\bm{b}=(b_{i})$ from $\mathbb{X}^{n}$, and a scalar $x\in\mathbb{X}$, vector addition and scalar multiplication follow the rules
$$
\{\bm{a}\oplus\bm{b}\}_{i}
=
a_{i}\oplus b_{i},
\qquad
\{x\bm{a}\}_{i}
=
xa_{i}.
$$ 

A vector with all entries equal to zero is called the zero vector and denoted by $\mathbb{0}$.

A vector is regular if it has no zero elements.

With the above operations, the set of vectors $\mathbb{X}^{n}$ forms a semimodule over an idempotent semifield.

\begin{figure}[ht]
\setlength{\unitlength}{1mm}
\begin{center}

\begin{picture}(50,40)

\put(0,5){\vector(1,0){50}}
\put(5,0){\vector(0,1){40}}

\put(5,5){\thicklines\vector(1,3){8.5}}
\put(13.5,30.5){\line(0,-1){27}}

\put(5,5){\thicklines\vector(2,1){34}}
\put(39,22){\line(-1,0){35}}

\put(5,5){\thicklines\vector(4,3){34}}
\put(39,30.5){\line(-1,0){35}}
\put(39,30.5){\line(0,-1){26.5}}

\put(0,0){$0$}

\put(13,0){$b_{1}$}
\put(38,0){$a_{1}$}

\put(0,22){$a_{2}$}
\put(0,30){$b_{2}$}

\put(13,33){$\bm{b}$}

\put(42,21){$\bm{a}$}

\put(35,33){$\bm{a}\oplus\bm{b}$}
\end{picture}
\hspace{10\unitlength}
\begin{picture}(40,45)

\put(5,5){\vector(1,0){35}}
\put(10,0){\vector(0,1){45}}

\put(10,5){\thicklines\vector(1,3){5}}
\put(15,20){\line(0,-1){16}}
\put(15,20){\line(-1,0){6}}

\put(10,5){\thicklines\vector(2,3){20}}
\put(30,35){\line(0,-1){31}}
\put(30,35){\line(-1,0){21}}

\put(5,10){\line(1,1){30}}

\put(5,0){$0$}
\put(13,23){$\bm{a}$}
\put(32,32){$x\bm{a}$}

\put(4,19){$a_{2}$}
\put(2,35){$xa_{2}$}

\put(14,0){$a_{1}$}
\put(28,0){$xa_{1}$}

\end{picture}

\end{center}
\caption{Vector addition (left) and scalar multiplication (right) in $\mathbb{R}_{\max,+}^{2}$.}\label{F-VAMS}
\end{figure}
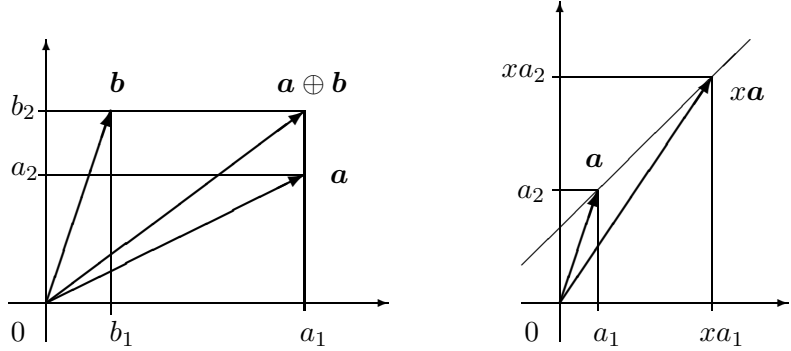

A geometric illustration for the operations in $\mathbb{R}_{\max,+}^{2}$ is given in Fig.~\ref{F-VAMS}.

Idempotent addition of two vectors in $\mathbb{R}_{\max,+}^{2}$ follows the ``rectangle rule'' that defines the sum as a diagonal of a rectangle formed by the coordinate axes together with lines drawn through the end points of the vectors. Scalar multiplication of a vector is equivalent to shifting the end point of the vector in the direction at $45^{\circ}$ to the axes.

As usual, a vector $\bm{b}\in\mathbb{X}^{n}$ is said to be linearly dependent on vectors $\bm{a}_{1},\ldots,\bm{a}_{m}\in\mathbb{X}^{n}$ if there are scalars $x_{1},\ldots,x_{m}\in\mathbb{X}$ such that
$$
\bm{b}
=
x_{1}\bm{a}_{1}\oplus\cdots\oplus x_{m}\bm{a}_{m}.
$$
In particular, $\bm{b}$ is collinear with $\bm{a}$ when $\bm{b}=x\bm{a}$.

Consider a set of vectors $\bm{a}_{1},\ldots,\bm{a}_{m}\in\mathbb{X}^{n}$. The set of all linear combinations
$$
\mathcal{A}
=
\{x_{1}\bm{a}_{1}\oplus\cdots\oplus x_{m}\bm{a}_{m}|x_{1},\ldots,x_{m}\in\mathbb{X}\}
$$
is referred to as the linear span of the vectors.

\begin{figure}[ht]
\setlength{\unitlength}{1mm}
\begin{center}

\begin{picture}(50,40)

\put(0,5){\vector(1,0){50}}
\put(5,0){\vector(0,1){40}}

\put(5,5){\thicklines\vector(1,3){4.25}}
\put(5,5){\thicklines\vector(2,3){17}}

\put(5,5){\thicklines\vector(4,1){23}}
\put(5,5){\thicklines\vector(2,1){34}}

\put(1.5,10){\thicklines\line(1,1){26}}
\multiput(2.5,11)(1,1){25}{\line(1,0){1}}

\put(17,0){\thicklines\line(1,1){26}}
\multiput(18,1)(1,1){25}{\line(-1,0){1}}

\put(39,30.5){\line(-1,0){35}}
\put(39,30.5){\line(0,-1){26.5}}

\put(10,18){\line(-1,0){6}}
\put(28,11){\line(0,-1){7}}

\put(5,5){\thicklines\vector(4,3){34}}

\put(0,0){$0$}

\put(7,21){$\bm{a}_{2}$}
\put(0,24){$x_{2}$}
\put(15,34){$x_{2}\bm{a}_{2}$}

\put(30,9){$\bm{a}_{1}$}
\put(32,0){$x_{1}$}
\put(42,20){$x_{1}\bm{a}_{1}$}

\put(31,34){$x_{1}\bm{a}_{1}\oplus x_{2}\bm{a}_{2}$}

\end{picture}

\end{center}
\caption{A linear span of two vectors in $\mathbb{R}_{\max,+}^{2}$.}\label{F-LC}

\end{figure}
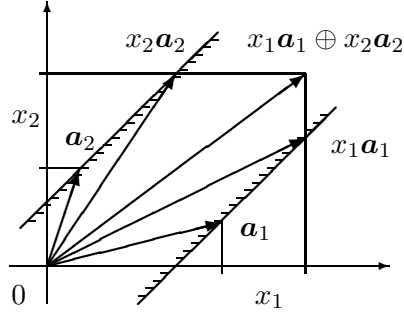

Specifically, the linear span of vectors $\bm{a}_{1}$ and $\bm{a}_{2}$ in $\mathbb{R}_{\max,+}^{2}$ has the form of a strip bounded by the lines drawn through the end points of the vectors (see Fig.~\ref{F-LC}).

For any column vector $\bm{x}=(x_{i})\in\mathbb{X}_{+}^{n}$, we introduce a row vector $\bm{x}^{-}=(x_{i}^{-})$ with elements $x_{i}^{-}=x_{i}^{-1}$ when $x_{i}\ne\mathbb{0}$, and $x_{i}^{-}=\mathbb{0}$ otherwise.

We define the distance between any two regular vectors $\bm{a}$ and $\bm{b}$ with a metric
$$
\rho(\bm{a},\bm{b})
=
\bm{b}^{-}\bm{a}\oplus\bm{a}^{-}\bm{b}.
$$

When $\bm{b}=\bm{a}$ we have $\rho(\bm{a},\bm{b})=\mathbb{1}$, where $\mathbb{1}$ is the minimum value the metric $\rho$ can take. 
 
Specifically, in $\mathbb{R}_{\max,+}^{n}$, we have $\mathbb{1}=0$, whereas the metric takes the form
$$
\rho(\bm{x},\bm{y})
=
\max_{1\leq i\leq n}|x_{i}-y_{i}|,
$$
and thus coincides with the classical Chebyshev metric.

For any conforming matrices $A=(a_{ij})$, $B=(b_{ij})$, and $C=(c_{ij})$ with entries in $\mathbb{X}$, matrix addition and multiplication together with multiplication by a scalar $x\in\mathbb{X}$ are performed in accordance with the formulas
\begin{gather*}
\{A\oplus B\}_{ij}
=
a_{ij}\oplus b_{ij},
\qquad
\{B C\}_{ij}
=
\bigoplus_{k}b_{ik}c_{kj},
\\
\{xA\}_{ij}=xa_{ij}.
\end{gather*}

A matrix with all entries equal to zero is called the zero matrix and denoted by $\mathbb{0}$.

A matrix is regular if it has no zero rows.

Consider the set of square matrices $\mathbb{X}^{n\times n}$. A matrix is diagonal if its off-diagonal entries are zero. The diagonal matrix $I=\mathop\mathrm{diag}(\mathbb{1},\ldots,\mathbb{1})$ is the identity matrix.

A matrix is reducible if it can be put in a block triangular form by simultaneous permutations of rows and columns. Otherwise, the matrix is irreducible.

For any matrix $A\ne\mathbb{0}$ and integer $p>0$, we have
$$
A^{0}
=
I,
\qquad
A^{p}
=
A^{p-1}A=AA^{p-1}.
$$

The trace of a matrix $A=(a_{ij})$ is defined as
$$
\mathop\mathrm{tr}A
=
\bigoplus_{i=1}^{n}a_{ii}.
$$

\subsection{Linear Operators and Linear Equations}

Any matrix $A\in\mathbb{X}^{m\times n}$ defines a mapping from the semimodule $\mathbb{X}^{n}$ to the semimodule $\mathbb{X}^{m}$. Since for any vectors $\bm{x},\bm{y}\in\mathbb{X}^{n}$ and scalar $\alpha\in\mathbb{X}$, it holds that
$$
A(\bm{x}\oplus\bm{y})
=
A\bm{x}\oplus A\bm{y},
\qquad
A(\alpha\bm{x})
=
\alpha A\bm{x},
$$
the mapping possesses the property of linear operators.

Suppose $A,C\in\mathbb{X}^{m\times n}$ are given matrices, and $\bm{b},\bm{d}\in\mathbb{X}^{m}$ are given vectors. A general linear equation in the unknown vector $\bm{x}\in\mathbb{X}^{n}$ is written in the form
$$
A\bm{x}\oplus\bm{b}
=
C\bm{x}\oplus\bm{d}.
$$ 

Note that due to the lack of additive inverse, one cannot put the equation in the form where all terms involving the unknown $\bm{x}$ are brought to one side of the equation while those without $\bm{x}$ go to another side.

Many practical problems reduce to solution of the following particular cases of the general equation
$$
A\bm{x}
=
\bm{d},
\qquad
A\bm{x}\oplus\bm{b}
=
\bm{x}.
$$

By analogy with linear integral equations, the above two equations are respectively referred to as that of the first kind and that of the second kind. The second-kind equations
$$
A\bm{x}
=
\bm{x}
\qquad
A\bm{x}\oplus\bm{b}
=
\bm{x}
$$
are also known in the literature as homogeneous and nonhomogeneous Bellman equations.

Some actual problems involve solution of inequalities of the first and second kinds in the form
$$
A\bm{x}
\leq
\bm{d}
\qquad
A\bm{x}\oplus\bm{b}\leq\bm{x}.
$$

\section{Preliminary Results}

Now we outline some recent results from \cite{Krivulin06Solution,Krivulin06Eigenvalues,Krivulin09Idempotent,Krivulin11Methods} that underlie subsequent applications of idempotent algebra to project scheduling problems.

\subsection{The First-Kind Equation and Inequality}

Given a matrix $A\in\mathbb{X}^{m\times n}$ and a vector $\bm{d}\in\mathbb{X}^{m}$, the problem is to find all solutions $\bm{x}\in\mathbb{X}^{n}$ of the equation and inequality given by
\begin{align}
A\bm{x}
&=
\bm{d},
\label{E-Axd}
\\
A\bm{x}
&\leq
\bm{d}.
\label{I-Axd}
\end{align}

A solution $\bm{x}_{0}$ to equation \eqref{E-Axd} is maximal if $\bm{x}_{0}\geq\bm{x}$ for all solutions $\bm{x}$ of \eqref{E-Axd}.

We present solution of equation \eqref{E-Axd} based on the analysis of the distance between vectors in $\mathbb{X}^{m}$. The solution involves the introduction of a new symbol
$$
\Delta
=
(A(\bm{d}^{-}A)^{-})^{-}\bm{d}
$$
to represent a residual quantity associated with \eqref{E-Axd}.

We start with a result that gives the distance from the vector $\bm{d}$ to the linear span of columns in the matrix $A$
$$
\mathcal{A}
=
\{A\bm{x}|\bm{x}\in\mathbb{X}^{n}\}.
$$

\begin{lemma}\label{L-minAxd}
Suppose $A\in\mathbb{X}^{m\times n}$ and $\bm{d}\in\mathbb{X}^{m}$ are regular matrix and vector. Then it holds that
$$
\min_{\bm{x}\in\mathbb{X}^{n}}\rho(A\bm{x},\bm{d})=\Delta^{1/2}
$$
with the minimum attained at
$$
\bm{x}
=
\Delta^{1/2}(\bm{d}^{-}A)^{-}.
$$
\end{lemma}

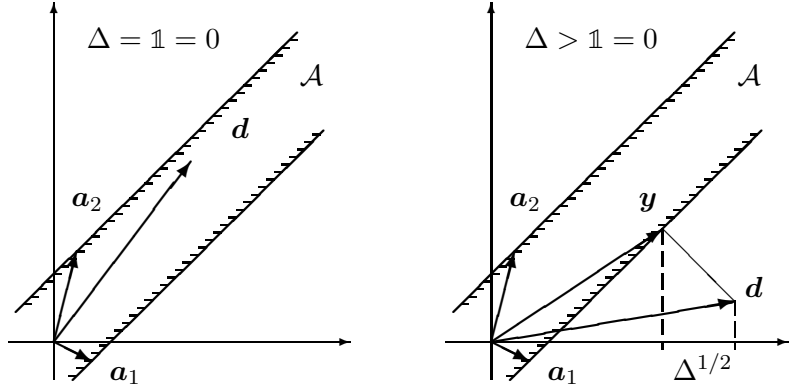
\begin{figure}[ht]
\setlength{\unitlength}{1mm}
\begin{center}

\begin{picture}(45,50)

\put(0,5){\vector(1,0){45}}
\put(6,0){\vector(0,1){50}}

\put(6,5){\thicklines\vector(1,4){3}}

\put(6,5){\thicklines\vector(2,-1){5}}

\put(1,9){\line(1,1){37}}
\multiput(2,10)(1,1){36}{\line(1,0){1}}
\put(1,9){\thicklines\line(1,1){37}}

\put(8.5,0){\line(1,1){33}}
\multiput(9,0.5)(1,1){33}{\line(-1,0){1}}
\put(8.5,0){\thicklines\line(1,1){33}}

\put(6,5){\thicklines\vector(3,4){18}}

\put(13,0){$\bm{a}_{1}$}
\put(8,23){$\bm{a}_{2}$}
\put(29,32){$\bm{d}$}

\put(38,39){$\mathcal{A}$}

\put(10,44){$\Delta=\mathbb{1}=0$}

\end{picture}
\hspace{10\unitlength}
\begin{picture}(45,50)

\put(0,5){\vector(1,0){45}}
\put(6,0){\vector(0,1){50}}

\put(6,5){\thicklines\vector(1,4){3}}

\put(6,5){\thicklines\vector(2,-1){5}}

\put(1,9){\line(1,1){37}}
\multiput(2,10)(1,1){36}{\line(1,0){1}}
\put(1,9){\thicklines\line(1,1){37}}

\put(8.5,0){\line(1,1){33}}
\multiput(9,0.5)(1,1){33}{\line(-1,0){1}}
\put(8.5,0){\thicklines\line(1,1){33}}

\put(28.5,20){\line(1,-1){9.5}}

\put(6,5){\thicklines\vector(3,2){22.5}}

\put(6,5){\thicklines\line(6,1){32}}
\put(36,10){\thicklines\vector(4,1){2}}

\multiput(28.5,20)(0,-3){5}{\line(0,-1){2}}
\put(28.5,5){\line(0,-1){1}}

\multiput(38,10.5)(0,-2.9){2}{\line(0,-1){2}}
\put(38,5){\line(0,-1){1}}

\put(13,0){$\bm{a}_{1}$}
\put(8,23){$\bm{a}_{2}$}
\put(39,11){$\bm{d}$}

\put(38,39){$\mathcal{A}$}

\put(25,23){$\bm{y}$}

\put(10,44){$\Delta>\mathbb{1}=0$}

\put(29.5,0){$\Delta^{1/2}$}

\end{picture}

\end{center}
\caption{A linear span $\mathcal{A}$ and a vector $\bm{d}$ in $\mathbb{R}_{\max,+}^{2}$ when $\Delta=\mathbb{1}$ (left) and $\Delta>\mathbb{1}$ (right).}\label{F-Lb}
\end{figure}

Fig.~\ref{F-Lb} presents examples of mutual arrangement of a vector $\bm{d}$ and the linear span $\mathcal{A}$ of columns $\bm{a}_{1}$ and $\bm{a}_{2}$ in a matrix $A$. In the case when $\Delta>\mathbb{1}$, the minimum distance to $\mathcal{A}$ is attained at the vector $\bm{y}=\Delta^{1/2}A(\bm{d}^{-}A)^{-}$. 

Furthermore, we consider sets
\begin{align*}
\mathcal{A}_{1}
&=
\{A\bm{x}|A\bm{x}\leq\bm{d},\bm{x}\in\mathbb{X}^{n}\},
\\
\mathcal{A}_{2}
&=
\{A\bm{x}|A\bm{x}\geq\bm{d},\bm{x}\in\mathbb{X}^{n}\}.
\end{align*}

\begin{lemma}\label{L-minAxdminAxd}
Suppose $A\in\mathbb{X}^{m\times n}$ and $\bm{d}\in\mathbb{X}^{m}$ are regular matrix and vector. Then it holds that
$$
\min_{A\bm{x}\leq\bm{d}}\rho(A\bm{x},\bm{d})
=
\min_{A\bm{x}\geq\bm{d}}\rho(A\bm{x},\bm{d})
=
\Delta,
$$
where the minimums are respectively attained at
$$
\bm{x}_{1}
=
(\bm{d}^{-}A)^{-},
\qquad
\bm{x}_{2}
=
\Delta(\bm{d}^{-}A)^{-}.
$$
\end{lemma}

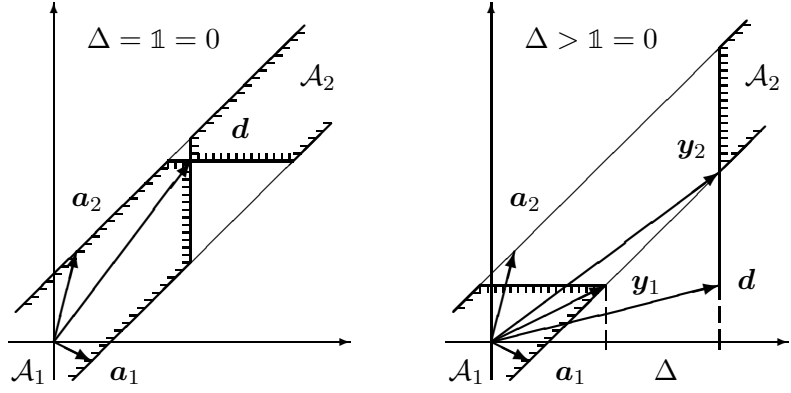
\begin{figure}[ht]
\setlength{\unitlength}{1mm}
\begin{center}

\begin{picture}(45,50)

\put(0,5){\vector(1,0){45}}
\put(6,0){\vector(0,1){50}}

\put(6,5){\thicklines\vector(1,4){3}}

\put(6,5){\thicklines\vector(2,-1){5}}

\put(1,9){\line(1,1){37}}
\multiput(2,10)(1,1){19}{\line(1,0){1}}
\put(1,9){\thicklines\line(1,1){20}}

\multiput(24,32)(1,1){15}{\line(1,0){1}}
\put(24,32){\thicklines\line(1,1){15}}

\put(8.5,0){\line(1,1){32}}
\multiput(9,0.5)(1,1){16}{\line(-1,0){1}}
\put(8.5,0){\thicklines\line(1,1){15.5}}

\multiput(38.5,30)(1,1){5}{\line(-1,0){1}}
\put(37.5,29){\thicklines\line(1,1){5}}

\put(6,5){\thicklines\vector(3,4){18}}

\put(21,29){\line(1,0){16.5}}
\put(21,29){\thicklines\line(1,0){16.5}}
\multiput(21,29)(1,0){3}{\line(0,-1){1}}
\multiput(24,29)(1,0){14}{\line(0,1){1}}

\put(24,15.5){\line(0,1){16.5}}
\put(24,15.5){\thicklines\line(0,1){16.5}}
\multiput(24,15.5)(0,1){14}{\line(-1,0){1}}
\multiput(24,30)(0,1){3}{\line(1,0){1}}

\put(13,0){$\bm{a}_{1}$}
\put(8,23){$\bm{a}_{2}$}
\put(29,32){$\bm{d}$}

\put(0,0){$\mathcal{A}_{1}$}
\put(38,39){$\mathcal{A}_{2}$}

\put(10,44){$\Delta=\mathbb{1}=0$}

\end{picture}
\hspace{10\unitlength}
\begin{picture}(45,50)

\put(0,5){\vector(1,0){45}}
\put(6,0){\vector(0,1){50}}

\put(6,5){\thicklines\vector(1,4){3}}

\put(6,5){\thicklines\vector(2,-1){5}}

\put(1,9){\line(1,1){38}}
\multiput(1.5,9.5)(1,1){3}{\line(1,0){1}}
\put(1,9){\thicklines\line(1,1){3.6}}

\multiput(36,44)(1,1){4}{\line(1,0){1}}
\put(36,44){\thicklines\line(1,1){3.6}}

\put(8.5,0){\line(1,1){32}}

\multiput(9,0.5)(1,1){12}{\line(-1,0){1}}
\put(8.5,0){\thicklines\line(1,1){12.5}}

\multiput(37.5,29)(1,1){5}{\line(-1,0){1}}
\put(36,27.5){\thicklines\line(1,1){6}}

\multiput(4.5,12.5)(1,0){16}{\line(0,-1){1}}
\put(4.5,12.5){\thicklines\line(1,0){16.5}}

\multiput(36,44)(0,-1){16}{\line(1,0){1}}
\put(36,44){\thicklines\line(0,-1){16.5}}

\put(6,5){\thicklines\vector(4,3){30}}

\put(6,5){\thicklines\vector(2,1){15}}

\put(6,5){\thicklines\vector(4,1){30}}

\put(36,12.5){\line(0,1){31.5}}

\multiput(21.1,12.5)(0,-3){3}{\line(0,-1){2}}
\put(21.1,5){\line(0,-1){1}}

\multiput(36,12.6)(0,-3){3}{\line(0,-1){2}}
\put(36,5){\line(0,-1){1}}

\put(14,0){$\bm{a}_{1}$}
\put(8,23){$\bm{a}_{2}$}
\put(38,12){$\bm{d}$}

\put(0,0){$\mathcal{A}_{1}$}
\put(39,39){$\mathcal{A}_{2}$}

\put(24,12){$\bm{y}_{1}$}
\put(30,30){$\bm{y}_{2}$}

\put(10,44){$\Delta>\mathbb{1}=0$}

\put(27,0){$\Delta$}

\end{picture}

\end{center}
\caption{The sets $\mathcal{A}_{1}$ and $\mathcal{A}_{2}$, and the vector $\bm{d}$ in $\mathbb{R}_{\max,+}^{2}$ when $\Delta=\mathbb{1}$ (left) and $\Delta>\mathbb{1}$ (right).}\label{F-LLb}
\end{figure}

A geometric interpretation in $\mathbb{R}_{\max,+}^{2}$ is given in Fig.~\ref{F-LLb}.

Note that if $\Delta>\mathbb{1}$ then the minimum distance from $\bm{d}$ to $\mathcal{A}_{1}$ and $\mathcal{A}_{2}$ is attained at respective vectors
$$
\bm{y}_{1}
=
A(\bm{d}^{-}A)^{-},
\qquad
\bm{y}_{2}
=
\Delta A(\bm{d}^{-}A)^{-}.
$$

The next statement is a consequence of the above results.
\begin{theorem}\label{T-EAxb}
Suppose $A\in\mathbb{X}^{m\times n}$ and $\bm{d}\in\mathbb{X}^{m}$ are regular matrix and vector. Then the following statements hold.
\begin{description}
\item[(a)] A solution of equation \eqref{E-Axd} exists if and only if $\Delta=\mathbb{1}$.
\item[(b)] If solvable, the equation has the maximum solution
$$
\bm{x}
=
(\bm{d}^{-}A)^{-}.
$$
\end{description}
\end{theorem}

Suppose that $\Delta>\mathbb{1}$. In this case equation \eqref{E-Axd} has no solution. However, we can define a quasi-solution to \eqref{E-Axd} as a solution of the equation
$$
A\bm{x}=\Delta^{1/2}A(\bm{d}^{-}A)^{-},
$$
which is always exists and takes the form
$$
\bm{x}_{0}
=
\Delta^{1/2}(\bm{d}^{-}A)^{-}.
$$

The quasi-solution yields the minimum deviation between the vectors $\bm{y}=A\bm{x}$ and the vector $\bm{d}$ in the sense of the metric $\rho$. When $\Delta=\mathbb{1}$, the quasi-solution obviously coincides with the maximum solution.

Consider the problem of finding two vectors $\bm{x}_{1}$ and $\bm{x}_{2}$ that provide the minimum deviation between both sides of \eqref{E-Axd}, while satisfying the respective inequalities
$$
A\bm{x}\leq\bm{d},
\qquad
A\bm{x}\geq\bm{d}.
$$
A solution to the problem is readily given by Lemma~\ref{L-minAxdminAxd}.

Finally, we present the following statement.
\begin{lemma}
For any matrix $A\in\mathbb{X}^{m\times n}$ and vector $\bm{d}\in\mathbb{X}^{m}$, the solution to inequality \eqref{I-Axd} is given by
$$
\bm{x}
\leq
(\bm{d}^{-}A)^{-}.
$$
\end{lemma}

The general solution to equation \eqref{E-Axd} with arbitrary matrix $A$ and vector $d$ is considered in \cite{Krivulin09Idempotent,Krivulin11Methods}.

\subsection{Second-Kind Equations and Inequalities}

Suppose a matrix $A\in\mathbb{X}^{n\times n}$ and a vector $\bm{b}\in\mathbb{X}^{n}$ are given, whereas $\bm{x}\in\mathbb{X}^{n}$ is an unknown vector. We examine the equation and inequality that have the form
\begin{align}
A\bm{x}\oplus\bm{b}
&=
\bm{x},
\label{E-Axbx}
\\
A\bm{x}\oplus\bm{b}
&\leq
\bm{x}.
\label{I-Axbx}
\end{align}

To solve equation \eqref{E-Axbx} we propose an approach based on the use of a function $\mathop\mathrm{Tr}(A)$ that takes each square matrix $A$ to a scalar according to the definition
$$
\mathop\mathrm{Tr}(A)
=
\bigoplus_{m=1}^{n}\mathop\mathrm{tr} A^{m}.
$$

The function is exploited to examine whether the equation has a unique solution, many solutions, or no solution, and so may play the role of the determinant.

The solution involves evaluation of matrices $A^{\ast}$, $A^{\times}$, and $A^{+}$. The matrices $A^{\ast}$ and $A^{\times}$ are given by 
$$
A^{\ast}
=
I\oplus A\oplus\cdots\oplus A^{n-1},
\quad
A^{\times}
=
A\oplus\cdots\oplus A^{n}.
$$

Let $\bm{a}_{i}^{\times}$ be column $i$ in $A^{\times}$, and $a_{ii}^{\times}$ be its diagonal element, $i=1,\ldots,n$. To construct the matrix $A^{+}$ we take the set of columns $\bm{a}_{i}^{\times}$ such that $a_{ii}^{\times}=\mathbb{1}$, and then reduce it by removing those columns that are linearly dependent on others. Finally, the columns in the reduced set are put together to form a matrix $A^{+}$.

We start with the solution of the homogeneous equation and inequality in the form
\begin{align}
A\bm{x}
&=
\bm{x},
\label{E-Axx}
\\
A\bm{x}
&\leq
\bm{x}.
\label{I-Axx}
\end{align}

The general solutions to the problems in the case of irreducible matrices are given by the following results.
\begin{lemma}
Let $\bm{x}$ be the solution of equation \eqref{E-Axx} with an irreducible matrix $A$. Then the following statements hold.
\begin{description}
\item[(a)] If $\mathop\mathrm{Tr}(A)=\mathbb{1}$, then $\bm{x}=A^{+}\bm{v}$ for any vector $\bm{v}$.
\item[(b)] If $\mathop\mathrm{Tr}(A)\ne\mathbb{1}$, then there is only the trivial solution $\bm{x}=\bm{\mathbb{0}}$.
\end{description}
\end{lemma}

\begin{figure}[ht]
\setlength{\unitlength}{1mm}
\begin{center}
\begin{picture}(33,50)

\put(0,23){\vector(1,0){33}}
\put(14,0){\vector(0,1){50}}

\put(14,23){\thicklines\vector(-1,0){10}}
\put(1,20){\thicklines\line(1,1){25}}

\put(14,23){\thicklines\vector(-1,-4){5}}
\put(7,1){\line(1,1){25}}

\put(14,23){\thicklines\vector(1,2){10}}

\multiput(24,43)(0,-3.1){7}{\line(0,-1){2}}

\multiput(24,43)(-3,0){4}{\line(-1,0){2}}

\put(0,26){$\bm{a}_{2}$}
\put(3,4){$\bm{a}_{1}$}

\put(27,41){$\bm{x}$}

\put(20,19){$x_{1}$}

\put(8,43){$x_{2}$}

\put(10,19){$0$}

\end{picture}
\hspace{10\unitlength}
\begin{picture}(33,50)

\put(0,23){\vector(1,0){33}}
\put(16,0){\vector(0,1){50}}

\put(16,23){\thicklines\vector(-1,0){12}}
\put(1,20){\thicklines\line(1,1){25}}
\multiput(2,21)(1,1){24}{\line(1,0){1}}

\put(16,23){\thicklines\vector(0,-1){8}}
\put(7,6){\thicklines\line(1,1){24}}
\multiput(8,7)(1,1){23}{\line(-1,0){1}}

\put(16,23){\thicklines\vector(2,3){12}}

\multiput(28,41)(0,-2.8){7}{\line(0,-1){2}}

\multiput(28,41)(-2.8,0){5}{\line(-1,0){2}}

\put(0,26){$\bm{a}_{2}$}
\put(18,12){$\bm{a}_{1}$}

\put(29,42){$\bm{x}$}

\put(27,19){$x_{1}$}

\put(10,41){$x_{2}$}

\put(12,19){$0$}

\end{picture}
\end{center}
\caption{Examples of solutions for homogeneous equations in $\mathbb{R}_{\max,+}^{2}$.}\label{F-GSEAxx}
\end{figure}
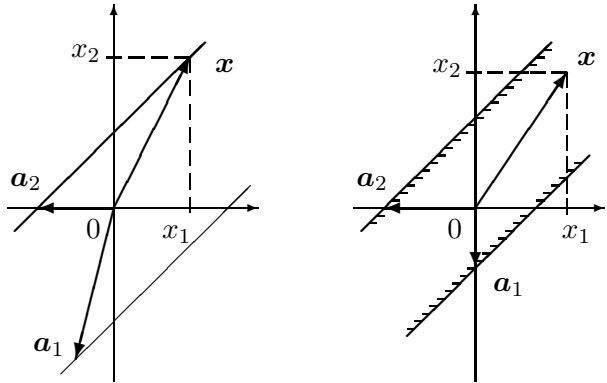

Fig.~\ref{F-GSEAxx} gives examples of solutions to homogeneous equations in $\mathbb{R}_{\max,+}^{2}$ for some particular matrices $A=(\bm{a}_{1},\bm{a}_{2})$ under the condition $\mathop\mathrm{Tr}(A)=\mathbb{1}$. On the left picture, the solution is shown with a thick line drawn through the end point of the vector $\bm{a}_{2}$. The right picture demonstrates the case when the solution coincides with the linear span of both columns in the matrix $A$.

\begin{lemma}
Let $\bm{x}$ be the solution of inequality \eqref{I-Axx} with an irreducible matrix $A$. Then the following statements hold.
\begin{description}
\item[(a)] If $\mathop\mathrm{Tr}(A)\leq\mathbb{1} $, then $\bm{x}=A^{\ast}\bm{v}$ for any vector $\bm{v}$.
\item[(b)] If $\mathop\mathrm{Tr}(A)>\mathbb{1}$, then there is only the trivial solution $\bm{x}=\bm{\mathbb{0}}$.
\end{description}
\end{lemma}

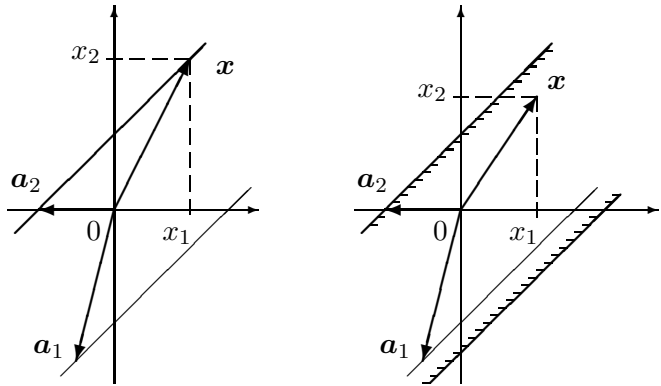
\begin{figure}[ht]
\setlength{\unitlength}{1mm}
\begin{center}
\begin{picture}(33,50)

\put(0,23){\vector(1,0){33}}
\put(14,0){\vector(0,1){50}}

\put(14,23){\thicklines\vector(-1,0){10}}
\put(1,20){\thicklines\line(1,1){25}}

\put(14,23){\thicklines\vector(-1,-4){5}}
\put(7,1){\line(1,1){25}}

\put(14,23){\thicklines\vector(1,2){10}}

\multiput(24,43)(0,-3.1){7}{\line(0,-1){2}}

\multiput(24,43)(-3,0){4}{\line(-1,0){2}}

\put(0,26){$\bm{a}_{2}$}
\put(3,4){$\bm{a}_{1}$}

\put(27,41){$\bm{x}$}

\put(20,19){$x_{1}$}

\put(8,43){$x_{2}$}

\put(10,19){$0$}

\end{picture}
\hspace{10\unitlength}
\begin{picture}(40,50)

\put(0,23){\vector(1,0){40}}
\put(14,0){\vector(0,1){50}}

\put(14,23){\thicklines\vector(-1,0){10}}
\put(1,20){\thicklines\line(1,1){25}}
\multiput(2,21)(1,1){24}{\line(1,0){1}}

\put(14,23){\thicklines\vector(-1,-4){5}}
\put(7,1){\line(1,1){25}}

\put(10,0){\thicklines\line(1,1){25}}
\multiput(10,0)(1,1){26}{\line(-1,0){1}}

\put(14,23){\thicklines\vector(2,3){10}}

\multiput(24,38)(0,-2.8){6}{\line(0,-1){2}}

\multiput(24,38)(-3,0){4}{\line(-1,0){2}}

\put(0,26){$\bm{a}_{2}$}
\put(3,4){$\bm{a}_{1}$}

\put(25,39){$\bm{x}$}

\put(20,19){$x_{1}$}

\put(8,38){$x_{2}$}

\put(10,19){$0$}

\end{picture}
\end{center}
\caption{Examples of solutions for a homogeneous equation (left) and inequality (right) in $\mathbb{R}_{\max,+}^{2}$.}\label{F-GSEIAxx}
\end{figure}

Fig~\ref{F-GSEIAxx} demonstrates solutions of homogeneous equation \eqref{E-Axx} and inequality \eqref{I-Axx} with a common matrix $A$.

In the general case of the nonhomogeneous equation and inequality, we have the following results.

\begin{theorem}\label{T-IMNHEGS}
Let $\bm{x}$ be the solution of equation \eqref{E-Axbx} with an irreducible matrix $A$. Then the following statements hold.
\begin{description}
\item[(a)] If $\mathop\mathrm{Tr}(A)<\mathbb{1}$, then $\bm{x}=A^{\ast}\bm{b}$.
\item[(b)] If $\mathop\mathrm{Tr}(A)=\mathbb{1}$, then $\bm{x}=A^{\ast}\bm{b}\oplus A^{+}\bm{v}$ for any vector $\bm{v}$.
\item[(c)] If $\mathop\mathrm{Tr}(A)>\mathbb{1}$, then $\bm{x}=\bm{\mathbb{0}}$ provided that $\bm{b}=\mathbb{0}$, and there is no solution otherwise.
\end{description}
\end{theorem}

\begin{lemma}
Let $\bm{x}$ be the solution of inequality \eqref{I-Axbx} with an irreducible matrix $A$. Then the following statements hold.
\begin{description}
\item[(a)] If $\mathop\mathrm{Tr}(A)\leq\mathbb{1} $, then $\bm{x}=A^{\ast}(\bm{b}\oplus\bm{v})$ for any vector $\bm{v}$.
\item[(b)] If $\mathop\mathrm{Tr}(A)>\mathbb{1}$, then $\bm{x}=\bm{\mathbb{0}}$ provided that $\bm{b}=\mathbb{0}$, and there is no solution otherwise.
\end{description}
\end{lemma}

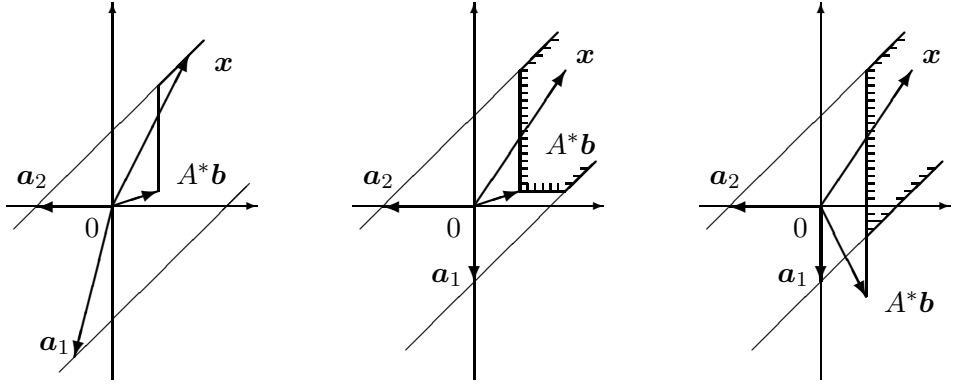
\begin{figure}[ht]
\setlength{\unitlength}{1mm}
\begin{center}

\begin{picture}(33,50)

\put(0,23){\vector(1,0){33}}
\put(14,0){\vector(0,1){50}}

\put(14,23){\thicklines\vector(-1,0){10}}
\put(20,39){\thicklines\line(1,1){6}}
\put(1,20){\line(1,1){25}}

\put(14,23){\thicklines\vector(-1,-4){5}}
\put(7,1){\line(1,1){25}}

\put(14,23){\thicklines\vector(1,2){10}}

\put(14,23){\thicklines\vector(3,1){6}}
\put(20,25){\thicklines\line(0,1){14}}

\put(1,26){$\bm{a}_{2}$}
\put(4,4){$\bm{a}_{1}$}

\put(27,41){$\bm{x}$}

\put(22,26){$A^{\ast}\bm{b}$}

\put(10,19){$0$}

\end{picture}
\hspace{10\unitlength}
\begin{picture}(33,50)

\put(0,23){\vector(1,0){33}}
\put(16,0){\vector(0,1){50}}

\put(16,23){\thicklines\vector(-1,0){12}}
\put(1,20){\line(1,1){25}}
\put(22,41){\thicklines\line(1,1){5}}
\multiput(22,41)(1,1){5}{\line(1,0){1}}

\put(16,23){\thicklines\vector(0,-1){10}}
\put(7,4){\line(1,1){25}}
\put(28,25){\thicklines\line(1,1){4}}
\multiput(29,26)(1,1){3}{\line(-1,0){1}}

\put(16,23){\thicklines\vector(2,3){12}}

\put(16,23){\thicklines\vector(3,1){6}}

\put(22,25){\thicklines\line(1,0){6}}
\multiput(22,25)(1,0){6}{\line(0,1){1}}

\put(22,25){\thicklines\line(0,1){16}}
\multiput(22,25)(0,1){16}{\line(1,0){1}}

\put(1,26){$\bm{a}_{2}$}
\put(10,13){$\bm{a}_{1}$}

\put(29,42){$\bm{x}$}

\put(25,30){$A^{\ast}\bm{b}$}

\put(12,19){$0$}

\end{picture}
\hspace{10\unitlength}
\begin{picture}(33,50)

\put(0,23){\vector(1,0){33}}
\put(16,0){\vector(0,1){50}}

\put(16,23){\thicklines\vector(-1,0){12}}
\put(1,20){\line(1,1){25}}
\put(22,41){\thicklines\line(1,1){5}}
\multiput(22,41)(1,1){5}{\line(1,0){1}}

\put(16,23){\thicklines\vector(0,-1){10}}
\put(7,4){\line(1,1){25}}
\put(22,19){\thicklines\line(1,1){10}}
\multiput(23,20)(1,1){9}{\line(-1,0){1}}

\put(16,23){\thicklines\vector(2,3){12}}

\put(16,23){\thicklines\vector(1,-2){6}}

\put(22,11){\thicklines\line(0,1){30}}

\multiput(22,20)(0,1){21}{\line(1,0){1}}

\put(1,26){$\bm{a}_{2}$}
\put(10,13){$\bm{a}_{1}$}

\put(29,42){$\bm{x}$}

\put(24,9){$A^{\ast}\bm{b}$}

\put(12,19){$0$}

\end{picture}

\end{center}
\caption{Examples of solutions for nonhomogeneous equations in $\mathbb{R}_{\max,+}^{2}$.}\label{F-GSEAxbx}
\end{figure}

A graphical illustration of the solution to the nonhomogeneous equations is given in Fig.~\ref{F-GSEAxbx}.

Related results for the case of arbitrary matrices can be found in \cite{Krivulin06Solution,Krivulin09Idempotent,Krivulin11Methods}.

\subsection{Eigenvalues and Eigenvectors}

A scalar $\lambda$ is an eigenvalue of a matrix $A\in\mathbb{X}^{n\times n}$ if there is a nonzero vector $\bm{x}\in\mathbb{X}^{n}$ such that
$$
A\bm{x}
=
\lambda\bm{x}.
$$

Any vector $\bm{x}\ne\mathbb{0}$ that satisfies the above equality is an eigenvector of $A$, corresponding to $\lambda$.

If the matrix $A\in\mathbb{X}^{n\times n}$ is irreducible, then it has only one eigenvalue given by
\begin{equation}
\lambda
=
\bigoplus_{m=1}^{n}\mathop\mathrm{tr}\nolimits^{1/m}(A^{m}).
\label{E-lambda}
\end{equation}

The corresponding eigenvectors of $A$ have no zero entries and take the form
$$
\bm{x}
=
A_{\lambda}^{+}\bm{v},
$$
where $A_{\lambda}=\lambda^{-1}A$ and $\bm{v}$ is any nonzero vector.

\begin{figure}[ht]
\setlength{\unitlength}{1mm}
\begin{center}

\begin{picture}(75,55)

\put(10,15){\vector(1,0){65}}
\put(20,5){\vector(0,1){50}}

\put(20,15){\thicklines\vector(4,1){20}}
\put(25,5){\line(1,1){42}}

\multiput(40,20)(-3,0){12}{\line(-1,0){2}}
\multiput(40,20)(0,-3.1){7}{\line(0,-1){2}}

\put(20,15){\thicklines\vector(3,2){45}}

\multiput(45,9)(0,-2.4){2}{\line(0,-1){1.9}}
\multiput(28.5,40)(-2.9,0){5}{\line(-1,0){2}}

\multiput(65,45)(0,-2.9){16}{\line(0,-1){2}}

\multiput(65,45)(-2.9,0){21}{\line(-1,0){2}}

\put(20,15){\thicklines\vector(1,3){8.3}}
\put(7,18.6){\thicklines\line(1,1){33}}

\put(20,15){\thicklines\vector(4,-1){25}}
\put(38.3,2){\thicklines\line(1,1){33}}

\put(16,45){\vector(0,-1){5}}
\put(16,44){\vector(0,1){1}}

\put(6,45){\vector(0,-1){25}}
\put(6,44){\vector(0,1){1}}

\put(46,6){\vector(1,0){19}}
\put(46,6){\vector(-1,0){1}}

\put(41,0){\vector(1,0){24}}
\put(41,0){\vector(-1,0){1}}

\put(30,35){$\bm{a}_{2}$}
\put(41,11.5){$\bm{a}_{1}$}

\put(38,22){$\bm{x}$}

\put(9,41){$x_{2}$}
\put(53,8){$x_{1}$}

\put(1,31){$\lambda$}

\put(51,1){$\lambda$}

\put(53,50){$A\bm{x}=\lambda\bm{x}$}

\end{picture}

\end{center}
\caption{An eigenvector $\bm{x}$ and eigenvalue $\lambda$ of a matrix $A$ in $\mathbb{R}_{\max,+}^{2}$.}\label{F-EV1}
\end{figure}
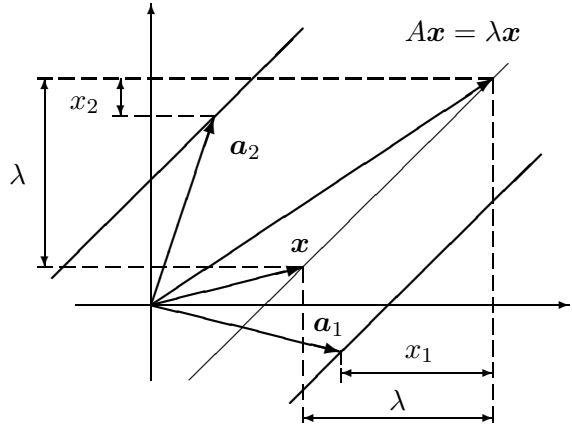

An example of an eigenvalue $\lambda$ and eigenvector $\bm{x}$ for a matrix $A=(\bm{a}_{1},\bm{a}_{2})$ in $\mathbb{R}_{\max,+}^{2}$ is given in Fig.~\ref{F-EV1}.

We conclude with an extremal property of the eigenvalue and eigenvectors of irreducible matrices. 

\begin{lemma}\label{L-mxmAx}
Suppose $A$ is an irreducible matrix with an eigenvalue $\lambda$. Then it holds that
$$
\min_{\bm{x}\in\mathbb{X}_{+}^{n}}\rho(A\bm{x},\bm{x})
=
\lambda\oplus\lambda^{-1}
$$
with the minimum attained at any eigenvector of $A$.
\end{lemma}

The eigenvalue-eigenvector problem and the above extremal property in the case of arbitrary matrices are examined in \cite{Krivulin06Eigenvalues,Krivulin09Idempotent,Krivulin11Methods}.

\section{Applications to Project Scheduling}

In this section we show how to apply results presented above to solve scheduling problems under various constraints (for further details on the schedule development in project management see, e.g., \cite{PMI08Aguide,Neumann03Project}).

As the underlying idempotent semiring, we use $\mathbb{R}_{\max,+}$ in all examples under discussion. 

\subsection{Start-to-Finish Precedence Constraints}

Consider a project that involves $n$ activities. Activity dependencies are assumed the form of Start-to-Finish relations that do not allow an activity to complete until some predefined time after initiation of other activities. The scheduling problem of interest consists in finding the latest initiation time for all activities subject to given constraints on their completion time.

For each activity $i=1,\ldots,n$, denote by $x_{i}$ its initiation time, and by $y_{i}$ its completion time. Let $d_{i}$ be a due date, and $a_{ij}$ be a minimum possible time lag between initiation of activity $j=1,\ldots,n$ and completion of $i$.

Given $a_{ij}$ and $d_{i}$, the completion time of activity $i$ must satisfy the relations
$$
y_{i}=d_{i},
\qquad
x_{j}+a_{ij}\leq y_{i},
\quad
j=1,\ldots,n.
$$ 

When $a_{ij}$ is not actually given for some $j$, it is assumed to be $\mathbb{0}=-\infty$.

The relations can be combined into one equation in the unknown variables $x_{1},\ldots,x_{n}$,
$$
\max(x_{1}+a_{i1},\ldots,x_{n}+a_{in})
=
d_{i}.
$$

By replacing the ordinary operations with those in $\mathbb{R}_{\max,+}$ in all equations, we get 
$$
a_{i1}x_{1}\oplus\cdots\oplus a_{in}x_{n}
=
d_{i},
\quad
i=1,\ldots,n.
$$

Furthermore, we introduce a matrix
$$
A
=
\left(
\begin{array}{ccc}
a_{11} & \ldots & a_{1n} \\
\vdots & \ddots & \vdots \\
a_{n1} & \ldots & a_{nn}
\end{array}
\right),
$$
and vectors
$$
\bm{d}
=
\left(
\begin{array}{c}
d_{1} \\
\vdots \\
d_{n}
\end{array}
\right),
\qquad
\bm{x}
=
\left(
\begin{array}{ccc}
x_{1} \\
\vdots \\
x_{n}
\end{array}
\right).
$$

The scheduling problem under the Start-to-Finish constraints leads us to the solution of the equation
$$
A\bm{x}
=
\bm{d}.
$$

Consider the residual $\Delta=(A(\bm{d}^{-}A)^{-})^{-}\bm{d}$ and suppose that $\Delta=\mathbb{1}=0$. According to Theorem~\ref{T-EAxb}, the equation has a maximum solution
$$
\bm{x}
=
(\bm{d}^{-}A)^{-}.
$$

If it appears that $\Delta>0$, then one can compute approximate solutions together with corresponding completion times as follows
\begin{align*}
\bm{x}_{0}
&=
\Delta^{1/2}(\bm{d}^{-}A)^{-},
&
\bm{y}_{0}
&=
A\bm{x}_{0};
\\
\bm{x}_{1}
&=
(\bm{d}^{-}A)^{-},
&
\bm{y}_{1}
&=
A\bm{x}_{1}
\leq
\bm{d};
\\
\bm{x}_{2}
&=
\Delta(\bm{d}^{-}A)^{-},
&
\bm{y}_{2}
&=
A\bm{x}_{2}
\geq
\bm{d}.
\end{align*}

Note that the completion times have their deviation from the due dates bounded with
$$
\rho(\bm{y}_{0},\bm{d})
=
\Delta^{1/2},
\qquad
\rho(\bm{y}_{1},\bm{d})
=
\rho(\bm{y}_{2},\bm{d})=\Delta.
$$

Suppose that the due date constraints are not mandatory and may be adjusted to some extent. As a vector of new due dates, it is natural to take $\bm{d}^{\prime}$ such that $\bm{y}_{1}\leq\bm{d}^{\prime}\leq\bm{y}_{2}$. In this case, deviation of the new due dates from the original ones does not exceed $\Delta$. The minimum deviation $\Delta^{1/2}$ is achieved when $\bm{d}^{\prime}=\bm{y}_{0}$.

As an example, consider a project with a constraint matrix
$$
A
=
\left(
\begin{array}{cccc}
8 & 10 & \mathbb{0} & \mathbb{0} \\
\mathbb{0} & 5 & 4 & 8 \\
6 & 12 & 11 & 7 \\
\mathbb{0} & \mathbb{0} & \mathbb{0} & 12
\end{array}
\right),
$$
and two due date vectors given by
$$
\bm{d}_{1}
=
\left(
\begin{array}{c}
14 \\
11 \\
16 \\
15
\end{array}
\right),
\qquad
\bm{d}_{2}
=
\left(
\begin{array}{c}
15 \\
15 \\
15 \\
15
\end{array}
\right).
$$

\begin{figure}[ht]
\setlength{\unitlength}{.75mm}
\begin{center}

\begin{picture}(100,40)

\put(5,35){\thicklines\circle{10}}
\put(2,34){$x_{1}$}
\put(35,35){\thicklines\circle{10}}
\put(32,34){$x_{2}$}
\put(65,35){\thicklines\circle{10}}
\put(62,34){$x_{3}$}
\put(95,35){\thicklines\circle{10}}
\put(92,34){$x_{4}$}

\put(5,5){\thicklines\circle{10}}
\put(2,4){$y_{1}$}
\put(35,5){\thicklines\circle{10}}
\put(32,4){$y_{2}$}
\put(65,5){\thicklines\circle{10}}
\put(62,4){$y_{3}$}
\put(95,5){\thicklines\circle{10}}
\put(92,4){$y_{4}$}

\put(5,30){\vector(0,-1){20}}
\put(9.5,32.5){\vector(2,-1){51}}

\put(35,30){\vector(0,-1){20}}
\put(31.5,31.5){\vector(-1,-1){23}}
\put(38.5,31.5){\vector(1,-1){23}}

\put(65,30){\vector(0,-1){20}}
\put(61.5,31.5){\vector(-1,-1){23}}

\put(95,30){\vector(0,-1){20}}
\put(90.5,32.5){\vector(-2,-1){51}}
\put(91.5,31.5){\vector(-1,-1){23}}

\put(1,22){$8$}
\put(10,26){$6$}
\put(22,29){$10$}
\put(36,22){$5$}
\put(42,29){$12$}
\put(54,29){$4$}
\put(59,22){$11$}
\put(87,22){$7$}
\put(96,22){$12$}
\put(78,29){$8$}

\end{picture}
\end{center}
\caption{An activity network with Start-to-Finish precedence relations.}\label{F-NM1}
\end{figure}
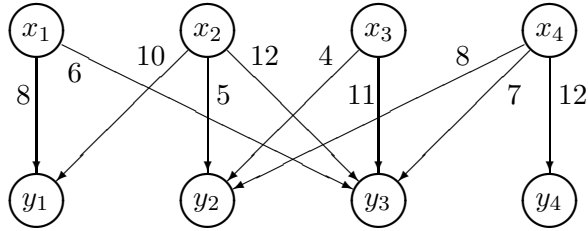

Fig.~\ref{F-NM1} demonstrates a network representation of the proceedings relations for activities in the project.

First we examine the equation $A\bm{x}=\bm{d}_{1}$. We have
$$
(\bm{d}_{1}^{-}A)^{-}
=
\left(
\begin{array}{c}
6 \\
4 \\
5 \\
3
\end{array}
\right),
\qquad
A(\bm{d}_{1}^{-}A)^{-}
=
\left(
\begin{array}{c}
14 \\
11 \\
16 \\
15
\end{array}
\right).
$$

Since $\Delta_{1}=(A(\bm{d}_{1}^{-}A)^{-})^{-}\bm{d}_{1}=0$, the equation has solutions including the maximum solution
$$
\bm{x}
=
(\bm{d}_{1}^{-}A)^{-}
=
(6, 4, 5, 3)^{T}.
$$

Now consider the equation $A\bm{x}=\bm{d}_{2}$. We get $\Delta_{2}=(A(\bm{d}_{2}^{-}A)^{-})^{-}\bm{d}_{2}=4>0$ and then conclude that the equation has no exact solutions. However, approximate solutions can be found as follows
\begin{align*}
&\bm{x}_{0}
=
(9,5,6,5)^{T},
\qquad
&&\bm{y}_{0}
=
(17,13,17,17)^{T},
\\
&\bm{x}_{1}
=
(7,3,4,3)^{T},
\qquad
&&\bm{y}_{1}
=
(15,11,15,15)^{T},
\\
&\bm{x}_{2}
=
(11,7,8,7)^{T},
\qquad
&&\bm{y}_{2}
=
(19,15,19,19)^{T}.
\end{align*}

\subsection{Start-to-Start Precedence Constraints}

Suppose there is a project consisting of $n$ activities and operating under Start-to-Start precedence constraints that determine the minimum allowed time intervals between initiation of activities. The problem is to find the earliest initiation time for each activity that does not violate these constraints.

For each activity $i=1,\ldots,n$, let $b_{i}$ be an early possible initiation time, and let $a_{ij}$ be a minimum possible time lag between initiation of activity $j=1,\ldots,n$ and initiation of $i$. The initiation time $x_{i}$ for activity $i$ is subject to the relations
$$
b_{i}\leq x_{i},
\qquad
a_{ij}+x_{j}\leq x_{i},
\quad
j=1,\ldots,n,
$$
where at least one must hold as an equality. 

We can replace the relations with one equation
$$
\max(x_{1}+a_{i1},\ldots,x_{n}+a_{in},b_{i})
=
x_{i}.
$$

Representation in terms of $\mathbb{R}_{\max,+}$, gives the scalar equations
$$
a_{i1}x_{1}\oplus\cdots\oplus a_{in}x_{n}\oplus b_{i}
=
x_{i},
\quad
i=1,\ldots,n.
$$

With the matrix-vector notation
$$
A
=
(a_{ij}),
\qquad
\bm{b}
=
(b_{1},\ldots,b_{n})^{T},
\qquad
\bm{x}
=
(x_{1},\ldots,x_{n})^{T}
$$
we arrive at a problem that is to solve the equation
$$
A\bm{x}\oplus\bm{b}
=
\bm{x}.
$$

Assume the matrix $A$ to be irreducible. It follows from Theorem~\ref{T-IMNHEGS} that if $\mathop\mathrm{Tr}(A)\leq\mathbb{1}=0$ then the equation has a solution given by
$$
\bm{x}
=
A^{\ast}\bm{b}\oplus A^{+}\bm{v},
$$
where $\bm{v}$ is any vector of appropriate size.

Consider a project with Start-to-Start relations and examine two cases, with and without early initiation time constraints. Let us define a matrix
$$
A
=
\left(
\begin{array}{rrrr}
0 & -2 & \mathbb{0} & \mathbb{0} \\
\mathbb{0} & 0 & 3 & -1 \\
-1 & \mathbb{0} & 0 & -4 \\
2 & \mathbb{0} & \mathbb{0} & 0
\end{array}
\right),
$$
and two vectors
$$
\bm{b}_{1}
=
\mathbb{0},
\qquad
\bm{b}_{2}
=
(1, 1, 2, 1)^{T}.
$$

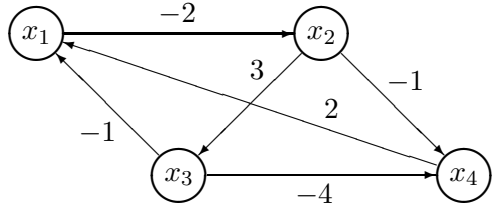
\begin{figure}[ht]
\setlength{\unitlength}{0.75mm}
\begin{center}

\begin{picture}(85,35)

\put(5,30){\thicklines\circle{10}}
\put(2,29){$x_{1}$}

\put(55,30){\thicklines\circle{10}}
\put(52,29){$x_{2}$}

\put(30,5){\thicklines\circle{10}}
\put(27,4){$x_{3}$}

\put(80,5){\thicklines\circle{10}}
\put(77,4){$x_{4}$}

\put(10,30){\vector(1,0){40}}

\put(26.5,8.5){\vector(-1,1){18}}

\put(58.5,26.5){\vector(1,-1){18}}

\put(51.5,26.5){\vector(-1,-1){18}}

\put(35,5){\vector(1,0){40}}

\put(75.25,6.7){\vector(-3,1){65.5}}

\put(12,11){$-1$}
\put(26,32){$-2$}
\put(42,22){$3$}
\put(50,0){$-4$}
\put(55,15){$2$}
\put(66,20){$-1$}

\end{picture}
\end{center}
\caption{An activity network with Start-to-Start precedence relations.}\label{F-NM2}
\end{figure}

A graph representation of the precedence relations involved in the project is depicted in Fig.~\ref{F-NM2}.

Let us first calculate the initiation time of activities in the project when $\bm{b}=\bm{b}_{1}=\mathbb{0}$ (that is, without early initiation time constraints imposed). Under this assumption, the equation becomes homogeneous and takes the form
$$
A\bm{x}
=
\bm{x}. 
$$

The matrix $A$ is irreducible and $\mathop\mathrm{Tr}(A)=0$. Therefore, the equation has a solution. 

Simple algebra gives
$$
A^{\ast}
=
A^{\times}
=
\left(
\begin{array}{rrrr}
0 & -2 & 1 & -3 \\
2 & 0 & 3 & -1 \\
-1 & -3 & 0 & -4 \\
2 & 0 & 3 & 0
\end{array}
\right).
$$

Note that all diagonal entries in $A^{\times}$ are equal to $\mathbb{1}=0$. However, considering that the first three columns are proportional, we take only one of them to form the matrix
$$
A^{+}
=
\left(
\begin{array}{rr}
-2 & -3 \\
0 & -1 \\
-3 & -4 \\
0 & 0
\end{array}
\right).
$$

The solution to the equation is given by
$$
\bm{x}
=
A^{+}\bm{v}
=
\left(
\begin{array}{rr}
-2 & -3 \\
0 & -1 \\
-3 & -4 \\
0 & 0
\end{array}
\right)
\bm{v},
\quad
\bm{v}\in\mathbb{R}_{\max,+}^{2}.
$$

Consider the case of the nonhomogeneous equation
$$
A\bm{x}\oplus\bm{b}_{2}
=
\bm{x}.
$$
We calculate the vector
$$
A^{\ast}\bm{b}_{2}
=
(3, 5, 2, 5)^{T},
$$
and then get
$$
\bm{x}
=
\left(
\begin{array}{c}
3 \\
5 \\
2 \\
5
\end{array}
\right)
\oplus
\left(
\begin{array}{rr}
-2 & -3 \\
0 & -1 \\
-3 & -4 \\
0 & 0
\end{array}
\right)\bm{v},
\quad
\bm{v}\in\mathbb{R}_{\max,+}^{2}.
$$

\subsection{Mixed Precedence Relations}

Consider a project that has both Start-to-Finish and Start-to-Start constraints. Let $A_{1}$ be a given Start-to-Finish constraint matrix, $\bm{d}$ a vector of due dates, and $\bm{x}$ an unknown vector of activity latest initiation time. To meet the constraints, the vector $\bm{x}$ must satisfy the inequality
$$
A_{1}\bm{x}
\leq
\bm{d}.
$$

Furthermore, there are also Start-to-Start constraints defined by a constraint matrix $A_{2}$. This leads to the equation
$$
A_{2}\bm{x}
=
\bm{x}.
$$

Suppose the equation has a solution $\bm{x}=A_{2}^{+}\bm{v}$. Substitution into the inequality yields
$$
A_{1}A_{2}^{+}\bm{v}
\leq
\bm{d}.
$$

Application of Theorem~\ref{T-EAxb} gives the maximum solution to the last inequality in the form $\bm{v}=(\bm{d}^{-}A_{1}A_{2}^{+})^{-}$. The solution to the whole problem is then written as
$$
\bm{x}
=
A_{2}^{+}(\bm{d}^{-}A_{1}A_{2}^{+})^{-}.
$$

As an illustration, we evaluate the solution to the problem under the condition that
$$
A_{1}
=
\left(
\begin{array}{cccc}
8 & 10 & \mathbb{0} & \mathbb{0} \\
\mathbb{0} & 5 & 4 & 8 \\
6 & 12 & 11 & 7 \\
\mathbb{0} & \mathbb{0} & \mathbb{0} & 12
\end{array}
\right),
\qquad
A_{2}
=
\left(
\begin{array}{rrrr}
0 & -2 & \mathbb{0} & \mathbb{0} \\
\mathbb{0} & 0 & 3 & -1 \\
-1 & \mathbb{0} & 0 & -4 \\
2 & \mathbb{0} & \mathbb{0} & 0
\end{array}
\right),
$$
and
$$
\bm{d}
=
(13,11,15,15)^{T}.
$$

By using results of previous examples, we successively get
$$
A_{1}A_{2}^{+}
=
\left(
\begin{array}{cc}
10 & 9 \\
8 & 8 \\
12 & 11 \\
12 & 12
\end{array}
\right),
\quad
(\bm{d}^{-}A_{1}A_{2}^{+})^{-}
=
\left(
\begin{array}{c}
3 \\
3
\end{array}
\right).
$$

Finally, we have
$$
\bm{x}
=
A_{2}^{+}(\bm{d}^{-}A_{1}A_{2}^{+})^{-}
=
(1, 3, 0, 3)^{T}.
$$

\subsection{Minimization of Maximum Flow Time}

Assume that a project has $n$ activities and operates under Start-to-Finish constraints. For each activity, consider the time interval between its initiation and completion, which is usually referred to as the flow time, the turnaround time or the processing time. The problem is to construct a schedule that minimizes the maximum flow time over all activities.

Let $A$ be an irreducible constraint matrix, $\bm{x}$ a vector of initiation time, and $\bm{y}=A\bm{x}$ a vector of completion time for the project. The problem can be formulated as that of finding a vector $\bm{x}$ that minimize
$$
\max(|y_{1}-x_{1}|,\ldots,|y_{n}-x_{n}|)
=
\rho(\bm{y},\bm{x}).
$$

In terms of $\mathbb{R}_{\max,+}$ we have
$$
\rho(\bm{y},\bm{x})
=
\rho(A\bm{x},\bm{x}).
$$

The problem of interest takes the form
$$
\min_{\bm{x}\in\mathbb{R}^{n}}\rho(A\bm{x},\bm{x})
$$
and can be solved by the application of Lemma~\ref{L-mxmAx}.

Let $\bm{d}$ be a given vector of activity due dates. Consider a problem of finding the latest initiation time for all activities so as to provide both the condition of minimum for the maximum flow time and the due date constraints in the form
$$
A\bm{x}
\leq
\bm{d}.
$$

By Lemma~\ref{L-mxmAx}, the first condition is satisfied when $\bm{x}$ is an eigenvector that corresponds to the eigenvalue $\lambda$ for the matrix $A$. The eigenvectors take the form $\bm{x}=A_{\lambda}^{+}\bm{v}$, where $A_{\lambda}=\lambda^{-1}A$ and $\bm{v}$ is any vector of appropriate size.  

By combining this result with the due date constraints, we get the inequality
$$
AA_{\lambda}^{+}\bm{v}
\leq
\bm{d}.
$$

With the maximum solution $\bm{v}=(\bm{d}^{-}AA_{\lambda}^{+})^{-}$ of the inequality, we arrive at the solution to the whole problem
$$
\bm{x}
=
A_{\lambda}^{+}(\bm{d}^{-}AA_{\lambda}^{+})^{-}.
$$

Let us evaluate the solution with the constraint matrix and due date vector defined as
$$
A
=
\left(
\begin{array}{ccc}
2 & 4 & 4 \\
2 & 3 & 5 \\
3 & 2 & 3
\end{array}
\right),
\qquad
\bm{d}
=
\left(
\begin{array}{c}
9 \\
8 \\
9
\end{array}
\right).
$$

First we get $\lambda=4$ with \eqref{E-lambda}, and define the matrix
$$
A_{\lambda}
=
\lambda^{-1}A
=
\left(
\begin{array}{rrr}
-2 & 0 & 0 \\
-2 & -1 & 1 \\
-1 & -2 & -1
\end{array}
\right).
$$

Furthermore, we have the matrices
$$
A_{\lambda}^{\ast}
=
A_{\lambda}^{\times}
=
\left(
\begin{array}{rrr}
0 & 0 & 1 \\
0 & 0 & 1 \\
-1 & -1 & 0
\end{array}
\right),
\quad
A_{\lambda}^{+}
=
\left(
\begin{array}{c}
1 \\
1 \\
0
\end{array}
\right),
$$
and then calculate
$$
AA_{\lambda}^{+}
=
\left(
\begin{array}{c}
5 \\
5 \\
4
\end{array}
\right),
\quad
(\bm{d}^{-}AA_{\lambda}^{+})^{-}
=
3.
$$

Finally, we arrive at the solution
$$
\bm{x}
=
A_{\lambda}^{+}(\bm{d}^{-}AA_{\lambda}^{+})^{-}
=
(4, 4, 3)^{T}.
$$

\section{Conclusion}
We have presented an approach that exploits idempotent algebra to solve computational problems in project scheduling. The approach allows to handle and combine different constraints and objectives that appear in actual problems in an easy and unified way. It is shown how to reformulate the problems in the algebraic setting, and then find related solutions based on recent results in the idempotent algebra theory. The solutions are given in a compact vector form that provides a basis for the development of efficient computation algorithms and software applications, including those intended for implementation on parallel and vector computers.

\bibliographystyle{utphys}

\bibliography{An_algebraic_approach_to_project_schedule_development_under_precedence_constraints}

\end{document}